\documentclass [12pt]{article}
\usepackage{amsmath,amssymb,amscd}
\usepackage{graphicx}

\begin{document}
\def\l{\lambda}
\def\tilde{\widetilde}
\def\m{\mu}
\def\a{\alpha}
\def\b{\beta}
\def\GCD{{\rm GCD }}
\def\g{\gamma}
\def\G{\Gamma}
\def\d{\delta}
\def\e{\epsilon}
\def\o{\omega}
\def\O{\Omega}
\def\v{\varphi}
\def\t{\theta}
\def\r{\rho}
\def\bs{$\blacksquare$}
\def\bp{\begin{proposition}}
\def\ep{\end{proposition}}
\def\bt{\begin{theo}}
\def\et{\end{theo}}
\def\be{\begin{equation}}
\def\ee{\end{equation}}
\def\bl{\begin{lemma}}
\def\el{\end{lemma}}
\def\bc{\begin{corollary}}
\def\ec{\end{corollary}}
\def\pr{\noindent{\bf Proof: }}
\def\note{\noindent{\bf Note. }}
\def\bd{\begin{definition}}
\def\ed{\end{definition}}
\def\C{{\mathbb C}}
\def\P{{\mathbb P}}
\def\Z{{\mathbb Z}}
\def\R{{\mathbb R}}
\def\d{{\rm d}}
\def\deg{{\rm deg\,}}
\def\deg{{\rm deg\,}}
\def\arg{{\rm arg\,}}
\def\min{{\rm min\,}}
\def\max{{\rm max\,}}
 \def\d{{\rm d}}
\def\deg{{\rm deg\,}}
\def\deg{{\rm deg\,}}
\def\arg{{\rm arg\,}}
\def\min{{\rm min\,}}
\def\max{{\rm max\,}}

\newtheorem{theo}{Theorem}[section]
\newtheorem{lemma}{Lemma}[section]
\newtheorem{definition}{Definition}[section]
\newtheorem{corollary}{Corollary}[section]
\newtheorem{proposition}{Proposition}[section]

\begin{titlepage}
\begin{center}

\topskip 2mm

{\LARGE{\bf {Parametric Center-Focus Problem

\medskip

for
Abel Equation }}}

\vskip 4mm

{\large {\bf M. Briskin $^{*}$,}}
\hspace {1 mm}
{\large {\bf F. Pakovich $^{**}$,}}
\hspace {1 mm}
{\large {\bf Y. Yomdin $^{***}$}}

\vspace{2 mm}
\end{center}

{$^{*}$ Jerusalem College of Engineering, Ramat Bet Hakerem, P.O.B. 3566,
Jerusalem 91035, Israel},
{$^{**}$ Department of Mathematics, Ben-Gurion University of the Negev,
Beer-Sheva 84105, Israel},
{$^{***}$ Department of Mathematics, The Weizmann Institute of
Science, Rehovot 76100, Israel}

\vspace{1 mm}

\noindent e-mail: briskin@jce.ac.il, \ pakovich@math.bgu.ac.il, \ yosef.yomdin@weizmann.ac.il

\vspace{1 mm}
\begin{center}

{ \bf Abstract}
\end{center}
{\small{

The Abel differential equation $y'=p(x)y^3 + q(x) y^2$ with meromorphic coefficients
$p,q$  is said to have a center on $[a,b]$ if all its solutions, with the
initial value $y(a)$ small enough, satisfy the condition $y(a)=y(b)$.
The problem of giving conditions on $(p,q,a,b)$ implying a center for the Abel equation
is analogous to the classical Poincar\'e Center-Focus problem for plane vector fields.

Following \cite{bfy1,bfy2,cgm1,cgm2} we say that Abel equation has a ``parametric center'' if for
each $\e \in \mathbb C$ the equation $y'=p(x)y^3 + \e q(x) y^2$ has a center. In the present paper
we use recent results of \cite{pak,bpy} to show show that for a {\it polynomial} Abel equation
parametric center implies strong ``composition'' restriction on $p$ and $q$. In particular, we show
that for $\deg p,q \leq 10$ parametric center is equivalent to the so-called ``Composition Condition''
(CC) (\cite{al,bfy1}) on $p,q$.

Second, we study trigonometric Abel equation, and provide a series of examples, generalizing a recent
remarkable example given in \cite{cgm1}, where certain moments of $p,q$ vanish while (CC) is violated.

\vspace{1 mm}
\begin{center}
------------------------------------------------
\vspace{1 mm}
\end{center}
This research was supported by the ISF, Grants No. 639/09 and  779/13, and by the
Minerva Foundation.

}}

\end{titlepage}
\newpage


\section{Introduction}
\setcounter{equation}{0}

We consider the Abel differential equation
\be \label{Abel}
y'=p(x)y^3 + q(x) y^2
\ee
with meromorphic coefficients $p,q$. A solution $y(x)$ of (\ref{Abel}) is called
``closed'' along a curve $\g$ with the endpoints $a,b$ if $y(a)=y(b)$ for the
initial element of $y(x)$ around $a$ analytically continued to $b$ along $\g$.
Equation (\ref{Abel}) is said to have a center along $\g$ if any its
solution $y(x)$, with the initial value $y(a)$ small enough, is closed
along $\g$. We always shall assume that $\g$ is a (closed or non-closed) curve
avoiding singularities of $p$ and $q$.
The Center-Focus problem for the Abel equation is to give a necessary and
sufficient condition on $p,q,\g$ for (\ref{Abel}) to have a center along $\g$. The
Smale-Pugh problem is to bound the number of isolated closed solutions of (\ref{Abel}).
The relation of these problems to the classical Hilbert 16-th and Poincar\'e
Center-Focus problems for plane vector fields is well known (see, eg. \cite{bfy1,broy} and references therein).

It turns out to be instructive to consider various parametric versions of the Center-Focus problem. In particular,
the following specific setting was considered in \cite{bfy2,bfy3,cgm1,cgm2}:

\bd \label{Par.Def}
Equation \ref{Abel} is said to have a ``parametric center'' if for each $\e \in \mathbb C$ the equation
\be \label{Abel.par}
y'=p(x)y^3 + \e q(x) y^2
\ee
has a center along $\g$.
\ed
It is easy to see that a parametric center can equivalently be defined by a requirement that the equation
$y'=\delta p(x)y^3 + q(x) y^2$ has a center for each $\delta$. 

The main purpose of this paper is to show that for a {\it polynomial} Abel
equation parametric center implies rather strong ``composition'' restriction on $p$ and $q$ (see Section 4). In particular, we show that for
$\deg p,q \leq 10$ parametric center is equivalent to ``Composition Condition'' (CC) on $p,q$ (see \cite{al,bfy1} and
Definition \ref{CC} below). For higher degrees $d$ of $p$ and $q$ we show that the dimension of possible ``non-composition'' 
couples $p,q$, forming parametric centers, is of order at most ${d}\over 3$, while the highest dimension of composition strata
is of order $d$.

Other parametric settings of the Center-Focus problem for Abel equation have been recently considered in \cite{ggl,cgm1,cgm2}
(see also \cite{pak4}). In particular, the following result has been obtained in \cite{cgm2}: if the equation
\be \label{Abel.par1}
y'=[\alpha p(x)+ \beta q(x)]y^3 + q(x) y^2
\ee
(with polynomial or trigonometric $p,q$) has a center for each $\alpha,\beta$ then $p,q$ satisfy the composition condition.
On the other hand, it was shown in \cite{cgm1} that for trigonometric polynomials $p,q$ 
the condition of
vanishing of the moments $\int P^kq, \int Q^kp$, closely related to the parametric center condition, 
does not imply the composition condition. In the last section of this paper, we construct a big series of examples, containing the examples given
in \cite{cgm1} as a particular case, where the vanishing of the moments $\int P^kq, \int Q^kp$ does not imply the composition
condition. 



\section{Poincar\'e mapping, Center Equations, and Composition condition}
\setcounter{equation}{0}

\subsection{Poincar\'e mapping and Center Equations} \label{Param.C.equat}

Both the Center-Focus and the Smale-Pugh  problems
can be naturally expressed in terms of the Poincar\'e
``first return'' mapping $y_b=G_\gamma (y_a)$ along $\gamma$. Let $y(x,y_a)$
denote the element around $a$ of the solution $y(x)$ of (\ref{Abel}) satisfying
$y(a)=y_a$. The Poincar\'e mapping $G_\gamma$ associates to each initial value
$y_a$ at $a$ the value $y_b$ at $b$ of the solution $y(x,y_a)$ analytically
continued along $\g$.

According to the definition above, the solution $y(x,y_a)$ is closed along
$\gamma$ if and only if $G_\gamma (y_a)=y_a$. Therefore closed solutions correspond
to the fixed points of $G_\gamma$, and (\ref{Abel}) has a center if and only if
$G_\gamma(y)\equiv y$. It is well known that $G_\gamma (y)$ for small $y$ is given
by a convergent power series

\be \label{eq1}
G_\gamma (y)= y + \sum_{k=2}^\infty v_k(p,q,\gamma)y^k.
\ee
Therefore the center condition $G_\gamma(y)\equiv y$ is equivalent to an
infinite sequence of algebraic equations on $p$ and $q$:

\be \label{eq2}
v_k(p,q,\gamma)=0, \ k=2,3,\dots.
\ee
Each $v_k(p,q,\gamma)$ can be expressed as a linear combination of certain iterated integrals along $\g$ of the form
$I_{\alpha}=\int h_{\alpha_1}\int h_{\alpha_2}\dots \int h_{\alpha_s}.$ Here $\alpha$ are the multi-indices
$\alpha =(\alpha_1, \dots , \alpha_s)$ with $\alpha_j = 1$ or $2$, and $h_1=p, \ h_2=q$. (see, for example, \cite{broy}).
The numbers $i,j$ of appearances of $p,q$ in these integrals satisfy $i+2j=k-1$, so $v_k(p,q,\gamma)$ are are weighted
homogeneous, but not homogeneous, polynomials in symbols $p,q$. The first few of $v_k(p,q,\gamma)$ are as follows:
\begin{eqnarray*}
v_2&=& -I_1, \ \ \ v_3 = 2I_{11} - I_2, \ \ \ v_4 = -6I_{111} +3I_{12} +2I_{21}\\
v_5&=& 24I_{1111} -12I_{112} -8I_{121} -6I_{211} +3I_{22}\\
v_6&=& -120I_{11111}+60I_{1112}+40I_{1121}+30I_{1211}+24I_{2111}\\
&-&15I_{122}-12I_{212}-8I_{221}\\
\end{eqnarray*}
Now, for parametric Abel equation (\ref{Abel.par}) we have:

\be \label{eq1.par}
G_\gamma (y,\e)= y + \sum_{k=2}^\infty v_k(p,q,\gamma,\e)y^k.
\ee
Consequently, we obtain the following result:

\bp\label{Par.C.eq}
The Center equations for parametric Abel equation (\ref{Abel.par}) take the following form:
\be \label{C.eq2.eps}
v_k(p,q,\gamma,\e)=\sum_{j=0}^{l(k)} v_{k,j}(p,q,\gamma)\e^j= 0, \ k=2,3,\dots,
\ee
where $v_{k,j}(p,q,\gamma)$ is a linear combination of iterated integrals as above, with exactly $j$ appearances of $q$,
and $l(k)=\lfloor {k\over 2} \rfloor-1$, where $\lfloor {k\over 2} \rfloor$ denotes the integer part of ${k\over 2}$.
\ep
\pr
This follows immediately from the description of the Center equations above. $\square$

\bc\label{Par.C.eq1}
Equation (\ref{Abel}) has a parametric center if and only if the following system of equations on $p,q,\g$ is satisfied:
\be\label{C.eq3.eps}
v_{k,j}(p,q,\gamma) = 0, \ k=2,3,\dots, \ j=0,1,\dots, l(k).
\ee
\ec
\pr
By definition, equation (\ref{Abel}) has a parametric center if and only if equation (\ref{Abel.par}) has a center for
each $\e$. By Proposition \ref{Par.C.eq} the polynomials $v_k(p,q,\gamma,\e)$ in $\e$ defined by equation (\ref{C.eq2.eps})
must vanish identically in $\e$. This is equivalent to vanishing of each of their coefficients $v_{k,j}(p,q,\gamma)$.
$\square$

\smallskip

So requiring center in (\ref{Abel.par}) to persist under variations of $\e$, we split the Center equations into their
summands corresponding to linear combinations of iterated integrals with the same numbers of the appearances of $q$.
While in general no description of these summands in ``closed form'' is known, the first two of them, $v_{k,1}, v_{k,2}$,
and the last one, $v_{k,l(k)},$ allow, at least partially, for such a description.

\medskip

To simplify notations, from now on, and till the end of Section \ref{Sec.Pol}, we always shall assume that $p,q$ are polynomials,
$\g=[a,b]$, and we shall denote by $P,Q$ the primitives $P(x)=\int_a^x p(\tau)d\tau$ and $Q(x)=\int_a^x q(\tau)d\tau$. Let
${\cal P}={\cal P}_{[a,b]}$ be the vector space of all complex polynomials $P$ satisfying $P(a)=P(b)=0$, and ${\cal P}_d$ the
subspace of $\cal P$ consisting of polynomials of degree at most $d$. We always shall assume that the polynomials

\be \label{iiuu}
P(x)=\int_a^x p(\tau)d\tau, \ \ \ Q(x)=\int_a^x q(\tau)d,
\ee
defined above are elements of $\cal P$. This restriction is natural in the study of the center conditions since it is
forced by the first two of the Center Equations.

\smallskip

The next result follows directly from the description of the Center equations given in \cite{broy} (compare also \cite{bpy},
Theorem 2.1).

\bt \label{ce.par}
I. The parametric Center equations $v_{k,j}(p,q,\gamma)=0$ for $j=0,1,2$ are as follows:

\smallskip

1. For $j=0$ and each $k$ the expressions $v_{k,0}(p,q,\gamma)$ vanish identically.

\smallskip

2. For $j=1$ and $k=0,1,2,$ we have $v_{k,1}\equiv 0.$ For $j=1$ and $k\geq 3,$
\be \label{mom}
v_{k,1}(P,Q)=m_{k-3}(P,Q)=\int_a^b P^{k-3}(x)q(x)dx=0.
\ee

3. For $j=2$ and $k\leq 5$ we have $v_{k,2}\equiv 0.$ For $j=2$ and $k\geq 6$ these equations are given by the coefficients of the
``second Melnikov function''
\be \label{mel}
v_{k,j}(P,Q)=D_k(P,Q)=0,
\ee
represented by linear combinations of iterated integrals in $p,q$ with exactly two appearances of $Q$ and $k-4$ appearances of
$P.$

\medskip

II. For $k\geq 2$, the highest index $j$ of the non-zero parametric Center equations is $j=l(k)=\lfloor {k\over 2} \rfloor-1.$

\smallskip

1. For $k\geq 2$ even, the corresponding equations are given by
\be \label{mominf}
v_{k,l(k)}(P,Q)=m_{l(k)}(Q,P)=\int_a^b Q^{l(k)}(x)p(x)dx=0.
\ee

2. For $k\geq 3$ odd, $v_{k,l(k)}$ are given by the coefficients of the ``second Melnikov function at infinity''
\be \label{melinf}
v_{k,l(k)}(P,Q)= \tilde D_k(P,Q)=0,
\ee
represented by linear combinations of iterated integrals in $p,q$ with exactly $l(k)$ appearances of $Q$ and two appearances of $P.$
\et
See Theorem \ref{Sec.Meln} in Section \ref{Sec.Pol} below, where some initial expressions $D_k(P,Q)$ are given explicitly (the
corresponding expressions $\tilde D_k(P,Q)$ ``at infinity'' are given in \cite{bpy}). Notice that considering the equation
$y'=\delta p(x)y^3 + q(x) y^2$ we get a slightly different set of parametric center equations.

As an immediate application of Theorem \ref{ce.par} we obtain the following known corollary (see \cite{bfy2}, \cite{bfy3}, \cite{cgm1}).

\bc\label{double}
Let Abel equation (\ref{Abel}) have a parametric center. Then the moment equations 
\be \label{doub}
m_i(P,Q)=\int_a^b P^i(x)q(x)dx=0, \ \  m_j(Q,P)=\int_a^b Q^j(x)p(x)dx=0,
\ee 
are satisfied for all $i,j\geq 0.$
\ec
\pr
Follows from Corollary \ref{Par.C.eq1} and Theorem \ref{ce.par} (cases I, 2, and II, 1).

\subsection{Center, Moment, and Composition Sets}

Let us consider $P,Q \in {\cal P}_d$. We define the parametric Center Set $PCS_d$ as the set of $(P,Q)\in {\cal P}_d\times {\cal P}_d$
for which equation (\ref{Abel}) has a parametric center. Equivalently, $PCS$ is the set of $(P,Q)$ satisfying Center Equations
(\ref{C.eq3.eps}). The moment set $MS^1_d$ (resp., $MS^2_d$)  consists of $(P,Q)\in {\cal P}_d\times {\cal P}_d$ satisfying Moment
equations (\ref{mom}) (resp., (\ref{mominf})). We put $MS_d=MS^1_d \cap MS^2_d$.

To introduce Composition Set $COS_d$ we recall the polynomial Composition Condition defined in \cite{bfy1}, which is a special case
of the general Composition Condition introduced in \cite{al} (for brevity below we will use the abbreviation ``CC'' for the
``Composition Condition'').

\bd \label{CC}
Polynomials $P,Q$ are said to satisfy the ``Composition Condition''  on $[a,b]$ if there exist polynomials $\tilde P$, $\tilde Q$ and
$W$ with $W(a)=W(b)$ such that $P$ and $Q$ are representable as $$P(x)=\tilde P(W(x)), \ \  \ Q(x)=\tilde Q(W(x)).$$ The Composition
Set $COS_d$ consists of all $(P,Q)\in {\cal P}_d\times {\cal P}_d$ satisfying the Composition Condition.
\ed
It is easy to see that the Composition Condition implies parametric center for (\ref{Abel}), as well as vanishing of each of the
moments and iterated integrals above. So we have $COS_d \subset PCS_d, \ COS_d \subset MS^1_d, \ COS_d \subset MS^2_d$.

\smallskip

It follows directly from Theorem \ref{ce.par} that the following statement is true:

\bp \label{inclusion}
$COS_d \subset PCS_d \subset MS_d=MS^1_d \cap MS^2_d.$
\ep
Our main goal will be to compare the parametric Center Set $PCS_d$ with the Composition Set $COS_d$. For this purpose we shall bound
the dimension of the non-composition components of $PCS_d$, analyzing the relations between $COS_d$ and the intersection $MS_d$.

\section{Moments vanishing and Composition}
\setcounter{equation}{0}

All the results in this section have been proved in \cite{bpy,pak}, so we state them below without proofs, and in a form convenient
for the purposes of the present paper.

\subsection{ $[a,b]$-Decompositions}
Let $a$ and $b$ be distinct points, and $P$ be a polynomial satisfying the condition $P(a)=P(b)$.
We are interested in ``$[a,b]$-decompositions'' of $P$, i.e. decompositions whose compositional ``right factors''
also take equal values at the points $a$ and $b$.

\bd \label{abcomp} Let a polynomial $P$ satisfying $P(a)=P(b)$ be given. We call a polynomial $W$ a right $[a,b]$-factor of $P$ if
$W(a)= W(b)$ and $P={\tilde P}\circ W$ for some polynomial $\tilde P$.
\ed
Recall that two decompositions  $P=P_1\circ W_1$ and $P=P_2\circ W_2$ of a polynomial $P$
into compositions of polynomials are called equivalent if there exists a polynomial $\mu$ of degree one such that
$$P_2=P_1\circ \mu^{-1}, \ \ \ W_2=\mu\circ W_1.$$ In accordance with this definition we shall call two right $[a,b]$-factors $W_1,$
$W_2$ of $P$ equivalent if $W_2=\mu\circ W_1$ for some polynomial $\mu$ of degree one.

\bd \label{abcomp+} A polynomial $P$ satisfying $P(a)=P(b)$ is called $[a,b]$-indecomposable if $P$ does not have right
$[a,b]$-factors non-equivalent to $P$ itself.
\ed

\noindent{\bf Remark.} Notice that any right $[a,b]$-factor of $P$ necessary has degree greater than one, and that
$[a,b]$-indecomposable $P$ may be decomposable in the usual sense.

\bp \label{iq} Any polynomial $P$ up to equivalence has a finite number of $[a,b]$-indecomposable right factors $W_j, \ \ j=1,\dots, s$.
Furthermore, each right $[a,b]$-factor $W$ of $P$ can be represented as $W=\tilde W(W_j)$ for some polynomial $\tilde W$ and
$j=1,\dots, s.$
\ep
It has been recently shown in \cite{pak} that for any polynomial $P$ the number $s$ of its non-equivalent $[a,b]$-indecomposable right
factors can be at most three. Moreover, if $s>1$ then these factors  necessarily have a very special form, similar to what appears in
Ritt's description in \cite{rit}.

The precise statement is given  by the following theorem (\cite{pak}, Theorem 5.3, \cite{bpy}, Theorem 3.1. Below
$T_d(x)=\cos(d \arccos x)$ denotes the $d$-th Chebyshev polynomial):

\bt \label{three} Let complex numbers $a\ne b$ be given. Then for any polynomial $P\in {\cal P}_{[a,b]}$ the number $s$ of its
$[a,b]$-indecomposable right factors $W_j$, up to equivalence, does not exceed 3.

Furthermore, if $s=2$, then either

$$P=U\circ z^{rn}R^n(z^n)\circ U_1, \ W_1=z^n\circ U_1,
\ W_2=z^r R(z^n)\circ U_1,$$ where
$R, U, U_1$ are polynomials, $r>0, n>1,$ $\GCD (n,r)=1,$
or

$$ P=U\circ T_{nm}\circ U_1, \ W_1=T_n\circ U_1, \ W_2=T_m\circ U_1,$$
where
$U, U_1$ are polynomials,
$n,m>1, \ \gcd (n,m)=1.$

On the other hand, if
$s=3$ then
$$P=U\circ z^2 R^2(z^2)\circ T_{m_1m_2}\circ U_1,$$ $$W_1=T_{2m_1}\circ U_1, \ \
\ W_2=T_{2m_2}\circ U_1, \ \ \  W_3=zR(z^2)\circ T_{m_1m_2}\circ U_1,$$
where
$R, U, U_1$ are polynomials,
$m_1,m_2 >1$
are odd, and $\GCD (m_1,m_2)=1.$

\et
In all the cases above we have $W_j(a)=W_j(b)$ and $U_1(a)\ne U_1(b)$.

\subsection{Moment vanishing versus $[a,b]$-Decompositions}

The main result of \cite{pak} can be formulated as follows:

\bt \label{yui} Let $P\in {\cal P}$  be given, and let $W_j,$ $j=1,\dots,s,$
be all its non-equivalent $[a,b]$-indecomposable right $[a,b]$-factors. Then
for any polynomial $Q$ all the moments $m_k=\int_a^b P^k(x)q(x)dx,$ $k\geq 0,$
vanish if and only if
$Q=\sum_{j=1}^s Q_j,$ where $Q_j=\tilde Q_j(W_j)$ for some polynomial $\tilde Q_j.$
\et
This theorem combined with Theorem \ref{three} provides an explicit description for
vanishing of the polynomial moments. In order to use it for the study of the Moment
Set, let us recall the notions of ``definite'' polynomials.

\bd \label{def.codef} A polynomial $P \in {\cal P}$ is called $[a,b]$-definite (or simply definite) if for any polynomial
$Q\in {\cal P}$ vanishing of the moments $m_k=\int_a^b P^k(x)q(x)dx,$ $k\geq 0,$ implies Composition Condition on $[a,b]$ for $P$ and $Q$.
\ed
Definite polynomials have been initially introduced and studied in \cite{pry}. Some their properties have been described in \cite{pp}.
The following theorems, providing a complete general characterization of definite polynomials, and their explicit description up to
degree $11$, have been obtained in \cite{bpy}.

\bt \label{definite} (\cite{bpy}) A polynomial $P$ is $[a,b]$-definite if and only if it has, up to equivalence, exactly one $[a,b]$-indecomposable
right factor $W$.
\et

Theorem \ref{definite} combined with Theorem \ref{three} allows us, at list in principle, to describe explicitly all the non-definite
polynomials up to a given degree. In particular, the following statement holds:

\bt \label{cla} (\cite{bpy})  For given $a\ne b$ non-definite polynomials $P \in {\cal P}_{11}$ appear only in degrees $6$ and $10$ and have, up to
change $P\rightarrow \lambda \circ P,$ where $\lambda$ is a polynomial of degree one, the following form:

\smallskip

1. $P_6=T_6 \circ \tau$, where $\tau$ is a polynomial of degree one transforming $a,b$ into $-{{\sqrt 3} \over 2}, {{\sqrt 3} \over 2}$.

\smallskip

2. $P_{10}=z^2 R^2(z^2)\circ \tau$, where $R(z)= z^2 + \gamma z + \delta$ is an arbitrary
quadratic polynomial satisfying $R(1)=0$ i.e. $\gamma+\delta=-1$, and $\tau$ is a
polynomial of degree one transforming $a,b$ into $-1,1$.
\et
The following definition introduces the zero subspace $Z(P)$ of the moments $m_i(P,Q)=\int_a^b P^i(x)q(x)dx$ for a given polynomial
$P$.

\bd \label{sums} Let $P \in {\cal P}$ be given. We define the set $Z(P)_d \subset {\cal P}_d$ as the set of polynomials
$Q \in {\cal P}_d$ which can be represented as $Q=\sum_{j=1}^s S_j(W_j),$ where $W_1,\dots, W_s$ are all $[a,b]$-indecomposable right
factors of $P$. The set $Z(P)$ is the union $\cup_d Z(P)_d$. Equivalently, $Z(P)_d$ consists of all $Q \in {\cal P}_d$ for which
$m_i(P,Q)=0, \ i=0,1,\ldots$.
\ed
For $P$ definite $Z(P)=\{S(W)\}$, where $W$ is the only $[a,b]$-indecomposable right factor of $P$. For a non-definite $P$ we have to consider sums
$Q=\sum_{j=1}^s S_j(W_j),$ where $W_1,\dots, W_s$ are all $[a,b]$-indecomposable right factors of $P$. The next result from \cite{bpy}
gives an example. To make formulas easier, here we shall assume that $[a,b]$ coincides with $[-{{\sqrt 3}\over 2},{{\sqrt 3}\over 2}]$.

\bt \label{noncodef} The set $Z(T_6)_d$ is a vector space consisting of all polynomials $Q\in {\cal P}_d$ representable as
$Q=S_1(T_2)+S_2(T_3)$ for some polynomials $S_1,$ $S_2$. The polynomials $S_1$ and $S_2$ can be chosen in such a way that $S_2$ is odd,
and $\max (2\,\deg S_1,3\,\deg S_2)\leq d$. Furthermore, the dimension ${\cal S}_{V,d}$ is equal to $[{{d+1}\over 2}]+[{{d+1}\over 3}]-[{{d+1}\over 6}].$
In particular, this dimension does not exceed $[{2\over 3}d]+1.$
\et

\section{Parametric Center for the Polynomial Abel Equation}\label{Sec.Pol}
\setcounter{equation}{0}

\bt\label{genpar}
Let Abel equation (\ref{Abel}) have a parametric center. Then either $P,Q$ satisfy Composition condition (CC) or both $P$ and $Q$ are
non-definite, $P\in Z(Q)$, and $Q\in Z(P)$.
\et
\pr
By Corollary \ref{double} if Abel equation (\ref{Abel}) has a parametric center, then
the moment equations \eqref{doub} are satisfied. Assume now
that $P,Q$ do not satisfy (CC). Then by Definition \ref{def.codef} both $P$ and $Q$ are non-definite, and by Definition \ref{sums} we have
$P\in Z(Q)$, and $Q\in Z(P)$. $\square$

\smallskip

Starting with \cite{bfy2,bfy3} the following conjecture has been discussed: {\it for a polynomial Abel equation parametric center is
equivalent to Composition condition (CC)}. A general description of parametric centers, provided by Theorem \ref{genpar} allows us to 
show that the set of $(P,Q)$ not satisfying (CC) but providing a parametric center to (\ref{Abel}) is ``small''. Below we give, based 
on Theorem \ref{genpar}, some specific results in this direction.

\bt\label{genpar1}
Let Abel equation (\ref{Abel}) with $P, \ Q \in {\cal P}_{11}$  have a parametric center. Then $P,Q$ satisfy Composition condition (CC).
\et
\pr
Assume that $P,Q$ do not satisfy (CC). Then by Theorem \ref{genpar} both $P$ and $Q$ are non-definite. Since the polynomials $P, \ Q$ are
contained in ${\cal P}_{11}$, they must have the form as provided by Theorem \ref{cla}. Although Theorem \ref{cla} provides two possibilities 
for each of the polynomials $P$ and $Q$, in any case $P$ has a right $[a,b]$-factor
$W_1=\alpha_1z^2+\beta_1z+\gamma_1$ of degree 2, and, similarly, $Q$ has a  right $[a,b]$-factor $W_2=\alpha_2z^2+\beta_2z+\gamma_2$ of degree 2.
Since the conditions $W_1(a)=W_1(b),$ $W_2(a)=W_2(b)$ yield the equalities 
$$\alpha_1(a+b)+\beta_1=0, \ \ \ \alpha_2(a+b)+\beta_2=0,$$ we conclude that 
$W_2=\lambda_1W_1+\lambda_2$ for some $\lambda_1,$ $\lambda_2 \in \C,$ and therefore both polynomials $P$ and $Q$ are polynomials in $W_1$, in contradiction with the initial condition.
$\square$

\bt\label{genpar2}
For $d\leq 11$ the non-composition components of the parametric Center set $PCS_d \subset {\cal P}_d\times {\cal P}_d$ is empty. For $d\geq 12$
the dimension of such components does not exceed $\lfloor {d\over 3} \rfloor+2$. In particular, this dimension is of order at most one third of the 
maximal dimension of the composition Center strata (which is of order $d$).
\et
\pr
The first part follows immediately from Theorem \ref{genpar1}. Let now assume that $d\geq 12$. By Theorem \ref{genpar}, for any non-composition couple
$(P,Q)$ both $P$ and $Q$ are non-definite. It was shown in \cite{bpy}, Proposition 3.3, that the dimension of the set of non-definite polynomials
in ${\cal P}_d$ does not exceed $\lfloor {d\over 6} \rfloor +1$. Hence the dimension of any non-composition components of $PCS_d$ does not exceed
$2(\lfloor {d\over 6} \rfloor +1)\leq \lfloor {d\over 3} \rfloor+2$. $\square$

\smallskip

\noindent{\bf Remark} One can hope to extend the approach of Theorem \ref{genpar1} to higher degrees of $P$ and $Q$.
Such an extension will require a better understanding of non-definite polynomials in higher degrees. In particular, the description of right
$[a,b]$-factors given by Theorem \ref{three} presents the following family of non-definite polynomials of degree $6m$: $P=T_6\circ U$, with $U$ an
arbitrary polynomial of degree $m$ transforming $a,b$ into ${{-\sqrt 3}\over 2},{{\sqrt 3}\over 2}$. The right factors of each such $P$ are
$W_1= T_2\circ U$ and $W_2= T_3\circ U$. It would be instructive to show that for (typical?) $P,Q$ in the above family the moments $\int P^kq$ and
$\int Q^kp$ cannot vanish simultaneously.

\smallskip

Theorems \ref{genpar}, \ref{genpar1}, \ref{genpar2} use only a very small part of the parametric center equations: the first and the last nonzero
ones, for each $k$. An interesting question is what consequences can be drawn from other combinations of the parametric center equations (see Theorem
\ref{ce.par}). The first (moment) equation for each $k$ corresponds to the ``infinitesimal center problem'', and its study was one of the main goals
of \cite{bfy1}-\cite{ggl}, \cite{mp1}-\cite{ppz}, and of many other publications. Recently in \cite{bpy} we've started the study of some situations
where we add more equations, in particular, those of the ``second Melnikov function at infinity'' in Theorem \ref{ce.par}. The approach of \cite{bpy}
can be used also in the study of the initial parametric center equations. Let us give the following definition:

\bd \label{Inf.l.Def}
Equation \ref{Abel} is said to have an ``infinitesimal center of order $l$'' if the Poincar\'e function $G_\gamma (y,\e)$ of the equation
$y'=p(x)y^3 + \e q(x) y^2$ satisfies ${{\partial^j}\over {\partial \e^j}}G_\gamma (y,\e)\equiv 0, \ j=0,\ldots,l-1$. Equivalently,
$v_{k,j}(p,q,\gamma)=0, \ k=2,3,\ldots, j=0,\ldots,l-1$.
\ed
\bt \label{Sec.Meln}
Let polynomial Abel equation \ref{Abel} have an infinitesimal center of order $3$. If $\deg p\leq 8, \deg q \leq 7,$ then $p,q$ satisfy Composition
condition (CC).
\et
\pr
By Definition \ref{Inf.l.Def} all the parametric center equations $v_{k,j}(p,q)=0$ are satisfied for $j=1,2$ and all $k$. Via Theorem \ref{ce.par}
those are the moments $m_{k-3}(P,Q)=\int_a^b P^{k-3}(x)q(x)dx, \ k=3,4,\ldots,$ and the coefficients of the second Melnikov function $D_k(P,Q)$.
These coefficients vanish identically for $k\leq 5$, and we shall use only three of them, namely (according to \cite{broy}):

\begin{eqnarray*}
D_6(P,Q) &=& {1\over 2}{\int_a^b pQ^2} \\
D_7(P,Q) &=& - {2\int_a^b PpQ^2} \\
D_8(P,Q) &=& \int_a^b P^3Qq - \\
-&&320\int_a^b P^2(t)q(t)dt\int_a^tPq+185\int_a^bP(t)q(t)dt
\int_a^t P^2q
\end{eqnarray*}
Now, let us assume that $P,Q$ do not satisfy (CC). Since the moment equations $m_j(P,Q)=0$ are satisfied, we conclude that $P$ is non-definite. The
only non-definite polynomial $P$ of degree less than $10$ is $T_6$. Hence $P=T_6\circ \tau$ where $\tau$ is a linear polynomial transforming $a,b$
into ${{-\sqrt 3}\over 2},{{\sqrt 3}\over 2}$. We conclude also that $Q$ belongs to the zero subspace $Z(P)$ (see Definition \ref{sums}). The
structure of $Z(P)$ is given by Theorem \ref{noncodef} above. Now we verbally repeat the arguments in the proof of Theorem 6.4 in \cite{bpy} (taking
into account that the coefficients of the iterated integrals in the expressions for $D_k(P,Q), \ k=6,7,8,$ given above, differ from the coefficients
in $\tilde D_k(P,Q)$ used in \cite{bpy}). $\square$

\smallskip

We conclude with another special example.
Let ${\cal R}=\{r_1,r_2,\dots\}$ be a set of prime numbers, finite or infinite.
Define $U({\cal R})$ as a subspace of ${\cal P}$ consisting of polynomials
$P=\sum_{i=0}^N a_i x^i$ such that for any non-zero coefficient $a_i$ the degree $i$
is either coprime with each $r_j\in {\cal R}$ or it is a power of some $r_j\in {\cal R}$.
Similarly, define $U_1({\cal R})$ as a subspace of ${\cal P}$ consisting of polynomials $Q$
such that for any non-zero coefficient $b_i$ of $Q$ all prime factors of $i$ are contained
in ${\cal R}$. In particular, if ${\cal R}$ coincides with the set of all primes numbers,
then $U({\cal R})$ consists of polynomials in ${\cal P}$ whose degrees with non-zero
coefficients are powers of primes, while $U_1({\cal R})={\cal P}.$

\bt\label{genpar4}
Let ${\cal R}=\{r_1,r_2,\dots\}$ be a set of prime numbers. Consider Abel equation (\ref{Abel}) with $P\in U({\cal R})$ and $Q\in U_1({\cal R}).$
Then this equation has a parametric center if and only if $P$ and $Q$ satisfy Composition condition (CC).
\et
\pr
It was shown in \cite{bpy} that all the polynomials $P\in U({\cal R})$ are definite with respect to polynomials $Q\in U_1({\cal R}).$ Hence the
result follows directly from Theorem \ref{genpar}. $\square$

\section{Trigonometric  Moments and Composition}\label{Sec.Trig}

Let $$P=P(\cos \theta,\sin \theta),\ \ \ Q=Q(\cos\theta ,\sin\theta )$$ now be trigonometric polynomials over $\R$, that is elements of the ring $\R_t[\theta]$ generated over $\R$ by the functions
$\cos \theta$, $\sin \theta$. Then in the same way as in the polynomial case the parametric center problem for Abel equation \eqref{Abel.par}
leads to the problem of characterization of $P,Q$ such that
\be \label{1} \int_0^{2\pi }P^idQ=0, \ \ \ i\geq 0,\ee
and
\be \label{1+} \int_0^{2\pi }Q^idP=0, \ \ \ i \geq 0.\ee
Indeed, our computations in Section \ref{Param.C.equat} above,  are essentially ``formal'', and can be applied to any required classes
of the coefficients $P,Q$, under minimal assumptions.

Again, a natural sufficient condition for \eqref{1} to be satisfied is related with compositional properties of $P$ and $Q$.
Namely, it is easy to see that if there exist $\widetilde P, \widetilde Q\in \R[x]$ and $W\in \R_t[\theta]$ such that
\be \label{c} P=\widetilde P\circ W, \ \ \ \ Q=\widetilde Q\circ W, \ee
then \eqref{1} holds. Furthermore, if for given $P$ there exist several such $Q$ (with different $W$),
then \eqref{1} obviously holds for their sum. In particular, the trigonometric moment problem \eqref{1} is closely related to the problem of description of solutions
of the equation
\be \label{mai} P=P_1\circ W_1=P_2\circ W_2,\ee where $P,W_1,W_2\in \R_t[\theta]$ and $P_1, P_2 \in \R[x],$
which was completely solved in the recent paper \cite{detri}.

Another sufficient condition for \eqref{1}, firstly proposed in \cite{ppre} in the complex setting, is the
following: the trigonometric polynomial $P$ has the form
\be \label{ur} P(\cos \theta,\sin\theta )=\widehat P(\cos (d\theta), \sin (d\theta))\ee
for some $\widehat P(x,y)\in \R[x,y]$  and $d>1$, while the trigonometric polynomial $Q$ is a linear combination of $\cos (m\theta)$ and $\sin (m\theta)$ for $m$ not divisible by $d$.
In this case, the sufficiency is a corollary  of the orthogonality of the functions $\cos k\theta,$ $\sin k\theta,$ $k\geq 1,$ on $[0,2\pi]$
(recall that any trigonometric polynomial
$P=P(\cos \theta,\sin \theta)$ can be represented as
$$P=\sum_{k=0}^r (a_k \cos (k\theta) + b_k \sin (k\theta))$$ for some $a_k,b_k \in \R$ and $r\geq 0$).

Notice that examples constructed in \cite{ppz}, \cite{abc} show that if we allow the coefficients of $P,Q$ to be complex numbers, then the two types of solutions  above (and their combinations) do not exhaust all possible solutions of \eqref{1}. On the other hand, for real $P,Q$ such examples are not known.

The second of the above mentioned sufficient conditions for \eqref{1} permits to construct examples of pairs $P,Q$ such that both equalities \eqref{1} and
\eqref{1+} hold but composition condition \eqref{c} does not hold. Indeed, let $d_1,d_2>1$ and $r_1,r_2\geq 1$ be integer numbers, and let
\be \label{fgh} P=\sum_{k=0}^{r_1} (a_k \cos (kd_1\theta) + b_k \sin (kd_1\theta)), \ \ \ \ Q=\sum_{l=0}^{r_2} (c_l \cos (ld_2\theta) + f_l \sin (ld_2\theta))\ee
be trigonometric polynomials  satisfying conditions $a_k=b_k=0$, whenever $d_2\vert k$, and
$c_l=f_l=0$, whenever $d_1\vert l$. Further, assume that $d_1$ and $d_2$ are coprime.
Then, since $P$ is a trigonometric polynomial in $\cos (d_1\theta)$, $\sin (d_1\theta)$
while $Q$ is some linear combination of $\cos (m\theta)$ and $\sin (m\theta)$ for $m$ not divisible by $d_1$,
equalities \eqref{1} hold. Similarly, equalities \eqref{1+} hold.

Further, it is not difficult to construct
pairs $P,Q$ as above for which condition \eqref{c} does not hold. Set for example $$P=\cos (d_1\theta), \ \  Q=\sin (d_2\theta).$$ It is easy to see that for any odd $d_1$ the derivatives of $P$ and $Q$, considered as functions of complex variable, have no common zeroes. On the other hand,
\eqref{c} implies that any complex zero of $W'$ is a common zero of $P'$ and $Q'$, and it is easy to see, using the formulas
$$\cos \theta=\frac{e^{i\theta}+e^{-i\theta}}{2}, \ \ \ \sin \theta=\frac{e^{i\theta}-e^{-i\theta}}{2i},$$
that any non-constant trigonometric polynomial $W$ with real coefficients has  complex zeroes (see e. g. \cite{pol}, Vol 2, Part 6). We conclude that 
\eqref{c} cannot be satisfied.


We can  modify the above construction
as follows.   Let $P=\cos (d_1\theta)$, while $Q$ be any trigonometric polynomial of the form \eqref{fgh}, where as above $c_l=f_l=0$, whenever $d_1\vert l$.
Set
$$\widetilde Q= Q+R\circ \cos (d_2\theta),$$ where $R$ is any polynomial of one variable. Then $\widetilde Q$ still has the form
$$\widetilde Q(\cos \theta,\sin\theta )=\widehat Q(\cos (d_2\theta), \sin (d_2\theta))$$
for some $\widehat Q(x,y)\in \R[x,y]$  implying that
$$\int_{0}^{2\pi}\widetilde Q^i d P=0,\ \ \ i\geq 0.$$ On the other hand,
equalities
$$\int_{0}^{2\pi}P^i d \widetilde Q=0,\ \ \ i\geq 0$$ hold by linearity, since for the polynomials
$$P=\cos (d_1\theta)=T_{d_1}\circ \cos \theta$$
and $$\widetilde Q- Q=R\circ \cos (d_2\theta)=R\circ
T_{d_2}\circ \cos \theta$$ condition \eqref{c} is satisfied for
$$P=T_{d_1}, \ \ \ Q=R\circ
T_{d_2}, \ \ \  W=\cos \theta.$$

The examples given in \cite{cgm1}, Proposition 18 are particular cases of the last series with
$$P=\cos (3\theta), \ \ \ \widetilde Q=\alpha\sin (2\theta)+\beta\cos (2\theta)+\gamma\cos (6\theta), \ \ \ \alpha, \beta, \gamma \in \R.$$
In this case, the easiest way to see that composition condition \eqref{c} does not hold is to observe that \eqref{c} implies the vanishing of  all integrals
\be \label{inte} \int_0^{2\pi}Q^id P^j, \ \ \ i\geq 0, j\geq 0,\ee while for given $P$ and $\widetilde Q$ integral \eqref{inte} is
distinct from zero for
$i=3,$ $j=2$, unless $a^2=3b^2$ (see  \cite{cgm1}).

\vskip1cm

\bibliographystyle{amsplain}

\end{document}